\documentclass[12pt]{article}
\usepackage{fancyhdr}
\usepackage{mathrsfs}
\usepackage{amsmath}
\usepackage{latexsym}
\usepackage{amssymb}
\usepackage{amsthm}
\usepackage{indentfirst}
\usepackage{titlesec}
\usepackage{titletoc}
\usepackage{graphicx}
\begin{document}
\title{ Classification of Ramification
Systems for Symmetric Groups}
\author{ Shouchuan Zhang $^a$,   \ \ Min Wu $^{a,  b}$,  \ \  Hengtai Wang $^{a}$ \\
a: Department  of Mathematics,  Hunan University,  Changsha,
410082\\b: Department of Mathematics,  Tsinghua University,
Beijing,  100084 }
\date{}


\newtheorem{Proposition}{\quad Proposition}[section]
\newtheorem{Theorem}[Proposition]{\quad Theorem}
\newtheorem{Definition}[Proposition]{\quad Definition}
\newtheorem{Corollary}[Proposition]{\quad Corollary}
\newtheorem{Lemma}[Proposition]{\quad Lemma}
\newtheorem{Example}[Proposition]{\quad Example}

\maketitle \addtocounter{section}{-1}

\numberwithin{equation}{section}

\begin {abstract}
 We  obtain the formula  computing the
number of isomorphic classes of ramification systems with
characters over group $S_n$ ($n \not=6$) and their
representatives.

\vskip 0.5cm
 \noindent 2000 Mathematics Subject Classification:
16W30,  05H99.

 \noindent Keywords:  Ramification,   character,   symmetry
group.
\end {abstract}

\section{\bf Introduction}

Classification of  Hopf algebras was developed and popularized in
the last decade of the twentieth century,   which would have
applications to a number of other areas of mathematics,  aside from
its intrinsic algebraic interest. In mathematical physics,
Drinfeld's and Jambo's work was to provide solutions to quantum
Yang-Baxter equation. In conformal field theory,  I. Frenkel and Y.
Zhu have shown how to assign a Hopf algebra to any conformal field
theory model\cite {FZ92}. In topology,  quasi-triangular and ribbon
Hopf algebras provide many invariants of knots,  links, tangles and
3-manifolds\cite {He91,  Ka97, Ra94, RT90}. In operator algebras,
Hopf algebras can be assigned as an invariant for certain
extensions.

Researches on the classification of Hopf algebra is in the
ascendant. The classification of monomial Hopf algebras,  which
are a class of co-path Hopf algebras,  and simple-pointed sub-Hopf
algebras of co-path Hopf algebras were recently obtained in \cite
{CHYZ04} and \cite {OZ04},  respectively. N. Andruskiewitsch and
H. J. Schneider have obtained interesting result in
classification of finite-dimensional pointed Hopf algebras
 with commutative coradical  \cite {AS98a,  AS98b,  AS02,  AS00}.
 More recently,  they  have also researched this problem in case of
 non-commutative coradical. Pavel Etingof and Shlomo Gelaki gave the complete
 and explicit classification of finite-dimensional
triangular Hopf algebras over an algebraically closed field $k$ of
characteristic $0$ \cite {EG00}.

Ramification systems with characters can be applied to classify PM
quiver Hopf algebras and multiple Taft algebras (see \cite[Theorem
3]{shou} ). In this paper,  we obtain the formula  computing the
number of isomorphic classes of ramification systems with
characters over group $S_n$ ($n \not=6$) and their
representatives.

\section{Preliminaries} \label {s2}

Unless specified otherwise,  in the paper we have the following
assumptions and notations. $n$ is a positive  integer with $n
\not=6$; $F$ is a field containing a primitive $n$-th root of 1;
$G= S_n$,  the symmetric group; Aut$G$ and Inn$G$ denote the
automorphism group and inner automorphism group,  respectively;
$S_I$ denotes the group of full permutations of set $I$; $Z(G)$
denotes the center of $G$ and $Z_x$ the centralizer of $x$; $1$
denotes the unity element of $G;$ $\widehat{G}$ denotes the set of
characters of all one-dimensional representations of $G$. $\phi_h$
denotes the inner automorphism induced by $h$ given by
$\phi_h(x)=hxh^{-1}, $ for any $x\in G$. $C_n$ is cyclic group of
order $n$; $D^B$ denotes the cartesian product $\prod \limits_
{i\in B}D_i$,  where $D_i = D$ for any $i \in B$.

\begin{Proposition}\label{7.2}
{(i)} There are
$\frac{n!}{\lambda_1!\lambda_2!\cdots\lambda_n!1^{\lambda_1}2^{\lambda_2}\cdots
n^{\lambda_n}}$ permutations of type
$1^{\lambda_1}2^{\lambda_2}\cdots n^{\lambda_n}$ in $S_n$;

{(ii)} If $n\neq 6$,   then Aut$(S_n)$=Inn$(S_n)$; If $n\neq 2,
6$,  then Inn$(S_n) \cong S_n$;

{(iii)} $S_n' = A_n$,  where $A_n$ is the alternating group. In
addition,
$$ S_n/A_n\cong\left\{
\begin{array}{lll}
C_2& n\geq2\\
C_1 &n =1
\end{array}
\right. . $$

\end{Proposition}

{\bf Proof. } (i) It follows from  \cite [ Exercise  2.7.7] {Hu}.

(ii) By  \cite [ Theorem 1.12.7] {Zh82},  $\mathrm{Aut}(A_n) \cong
S_n$ when $n >3$ and  $n \not=6$.  By  \cite [Theorem 1.6.10]
{Hu74},  $A_n $ is a non-commutative simple group when  $n \not=
4$ and $n>1$. Thus it follows from  \cite [Theorem 1.12.6 ] {Zh82}
that $\mathrm{Aut} (S_n) = \mathrm{Inn} (S_n)\cong S_n$ when  $n
\not= 6$ and $n>4$. Applying \cite [Example in Page 100] {Zh82},
we obtain
 $\mathrm{Aut} (S_n) = \mathrm{Inn} (S_n)\cong S_n$ when  $n=3,  4$.
Obviously,  $\mathrm{Aut} (S_n) = \mathrm{Inn} (S_n)$ when  $n=1,
2$.

 (iii) It follows from  \cite [Theorem 2.3.10] {Xu95}. $\Box$

\begin{Definition}\label{7.4}
A group $G$ is said  to act on a non-empty set $\Omega$,  if there
is a map  $G\times \Omega\rightarrow\Omega$,  denoted by $(g,
x)\mapsto g\circ x$,  such that for all $x\in \Omega$ and $g_1,
g_2\in G$:
$$
(g_1g_2)\circ x=g_1\circ(g_2\circ x)\qquad and\qquad   1\circ x=x,
$$ where 1 denote the unity element of $G.$

\end{Definition}

For each $x\in \Omega$,  let
$$ G_x:=\{g\circ x \mid  g\in G\}, $$
 called the orbit of $G$ on $\Omega$. For each $g\in G$,  let
 $$ F_g:=\{x\in \Omega \mid  g\circ x=x\}, $$
 called the  fixed point set of  $g$.

Burnside's lemma,  sometimes also called Burnside's counting
theorem,  which is  useful to compute the number of orbits.


\begin{Theorem}\label{7.7} (See \cite [Theorem 2.9.3.1]{Hu})
( Burnside's Lemma ) Let $G$ be a finite group that acts on  a
finite set $\Omega$,  then the number $\mathcal {N}$ of orbits is
given by the following formula:
$$
\mathcal {N}=\frac{1}{ \mid  G \mid  }\sum\limits_{g\in G} \mid F_g
\mid  . $$

\end{Theorem}


\begin{Definition}\label{7.8}

 Let $G$ be a  group. A character of a  one-dimensional representation of $G$ is  a homomorphism from $G$ to $F-\{0\}$,  that
is ,  a map $\chi:G\rightarrow F$ such that
$\chi(xy)=\chi(x)\chi(y)$ \ and\ $\chi(1)=1$ for any $ x, y\in G$.
The set of all characters of all one-dimensional representations
of $G$ is denoted by $\widehat{G}$. \end{Definition} For
convenience,  a character of a  one-dimensional representation of
$G$ is called a character of $G$ in short. \cite[Page 36, example
4]{char} indicated that there is a one-to-one correspondence
between characters  of $G$ and characters  of $G/G'$.

\begin{Proposition}\label{7.9}

 Let $G$ be a finite abelian group  and $F$ contains $\mid \! G \! \mid $th primitive  root of 1,  then
there are exactly $\mid \! G \! \mid $  characters of
 $G$,  i.e. $ \mid \widehat{G}
\mid  = \mid  G \mid  . $

\end{Proposition}

\begin{Corollary}\label{7.10}

 Let $G$ be a finite  group  and $F$ contains $\mid \! G \! \mid $th primitive  root of 1,  then
there are exactly $\mid \! G/G' \! \mid $  characters of $G$,
i.e. $ \mid \widehat{G} \mid  = \mid  G/G' \mid  . $

\end{Corollary}

\section{Classification of ramification systems for the symmetric group}
\subsection{Ramification Systems and Isomorphisms}

We begin with defining  ramification and ramification system.
\begin{Definition}\label{7.11}
Let $G$ be a finite  group and $\mathcal {K}(G)$ the set of all
 conjugate classes of $G$. $r$ is called a ramification data or ramification of $G$ if $r=\sum\limits_{C\in \mathcal
{K}(G)}r_{_C}C$,  where  $r_C$ is a non-negative integer for any
$C$ in $\mathcal {K}(G)$. For convenience,  we choice a set
$I_C(r)$ such that  $r_{_C}= \mid  I_C(r) \mid  $ for any  $C$ in
$\mathcal {K}(G)$.

Let $\mathcal {K}_r(G)=\{C\in \mathcal {K}(G) \mid  r_{_C}\neq
0\}=\{C\in \mathcal {K}(G) \mid  I_C(r)\neq \emptyset\}$.

\end{Definition}

\begin{Definition}\label{7.12}
$(G, r, \overrightarrow{\chi}, u)$ is called a ramification system
 with characters (or RSC in short), if r is a ramification of G, u is a
 map from $\mathcal{K}(G)$ to G with u(C)$\in$ C for any
 C$\in\mathcal{K}(G)$, and $\overrightarrow{\chi}=\{\chi_C^{(i)}\}_{i\in I_C(r), C\in
\mathcal {K}_r(G)}\in\prod\limits_{C\in \mathcal
{K}_r(G)}(\widehat{Z_{u(C)}})^{r_{_C}}$.

\end{Definition}

\begin{Definition}\label{7.13}

 $RSC(G, r, \overrightarrow{\chi}, u)$ and $
RSC(G', r', \overrightarrow{\chi'}, u')$ are said to be
isomorphic(denoted $RSC(G, r, \overrightarrow{\chi}, u)\cong
RSC(G', r', \overrightarrow{\chi'}, u')$) if the following
conditions are satisfied:

{(i)} There exists a group isomorphism $\phi:G\rightarrow G';$

{(ii)} For any $ C\in \mathcal {K}(G), $there exists an element
$h_C\in G$ such that $\phi(h_Cu(C)h_C^{-1})=u'(\phi(C))$;

{(iii)} For any $ C\in\mathcal {K}_r(G), $there exists a bijective
map $\phi_C:I_C(r)\rightarrow I_{\phi(C)}(r')$ such that
$\chi_{\phi(C)}^{\prime(\phi_C(i))}(\phi(h_Chh_C^{-1}))=\chi_C^{(i)}(h),
$ for all $h\in Z_{u(C)}, $ for all $i\in I_C(r)$.
\end {Definition}

{\bf Remark.} Obviously,  if $u(C)=u'(C)$ for any $ C\in \mathcal
{K}_r(G), $  then $RSC(G, r, \overrightarrow{\chi}, u)\cong RSC(G,
r, \overrightarrow{\chi}, u')$. In fact,  let $\phi=id_G$ and
$\phi_C=id_{I_C(r)}$ for any $ C\in \mathcal {K}_r(G)$. It is
straightforward to check the conditions in Definition\ref{7.13}
hold.

\begin{Lemma}\label{2.4}

Given an $RSC(G, r, \overrightarrow{\chi}, u)$.If $RSC(G', r',
\overrightarrow{\chi'}, u')\cong RSC(G, r, \overrightarrow{\chi},
u)$, then there exists a group isomorphism $\phi:G\rightarrow G'$
such that

{(i)} $ r'_{\phi(C)}=r_{_C}$ for any $ C\in \mathcal {K}(G)$;

{(ii)}  $\mathcal {K}_{r'}(G')=\phi(\mathcal {K}_{r}(G))=\{\phi(C)
\mid C\in \mathcal {K}_r(G)\}.$

\end{Lemma}
{\bf Proof.} (i) By Definition \ref{7.13},  there exists a group
isomorphism $\phi:G\rightarrow G'$ such that $r'_{\phi(C)}=r_{_C}$
for any  $ C\in \mathcal {K}_r(G)$. If $C\notin
 \mathcal {K}_r(G)$,  i.e. $r_{_C}=0, $ then $r'_{\phi(C)}=r_{_C}=0$;

(ii) is an immediate consequence of (i).
 $\Box$

\begin{Proposition}\label{2.5}
Given a $RSC(G, r, \overrightarrow{\chi}, u)$,  then $ RSC(G,  r',
\overrightarrow{\chi'},  u')\cong RSC(G,  r,
\overrightarrow{\chi},  u)$ if and only if there exists a group
isomorphism $\phi \in \hbox {Aut} G$ such that:

{(i)} $r'=\sum\limits_{C\in \mathcal {K}(G)}r_{\phi^{-1}(C)}C$;

{(ii)}
$\chi_{\phi(C)}^{\prime(i)}=\chi_C^{(\phi_C^{-1}(i))}(\phi\phi_{h_C})^{-1}$
and $u'(\phi(C))=\phi\phi_{h_C}(u(C))$ for any $ C\in \mathcal
{K}_r(G)$ ,  $i\in I_C(r) $,  $\phi\in \mathrm{Aut}G$ and
$\phi_C\in S_{r_{_C}}$,  $h_C\in G$.

\end{Proposition}
{\bf Proof.} If $RSC(G, r', \overrightarrow{\chi'}, u')\cong
RSC(G, r, \overrightarrow{\chi}, u) $,  then,  by Lemma
\ref{2.4}(i),  there exists $\phi\in \mathrm{Aut}G$ such that
$$r'=\sum\limits_{C\in
\mathcal
 {K}(G)}r'_CC=\sum\limits_{C\in\mathcal
 {K}(G)}r'_{\phi(C)}\phi(C)=\sum\limits_{C\in\mathcal
 {K}(G)}r_{_C}\phi(C)=\sum\limits_{C\in\mathcal {K}(G)}r_{\phi^{-1}(C)}C.
$$
This shows Part (i).

By Definition \ref{7.13} (ii),  there exists $h_C\in G$ such that
$u'(\phi(C))=\phi(h_Cu(C)h_C^{-1})=\phi\phi_{h_C}(u(C))$ for any $
C\in \mathcal {K}_r(G)$; By Definition \ref{7.13} (iii),  there
exists $\phi_C\in S_{I_C(r)}=S_{r_{_C}}$ such that
$\chi_{\phi(C)}^{\prime(\phi_C(i))}(\phi(h_Chh_C^{-1}))
=\chi_{\phi(C)}^{\prime(\phi_C(i))}\phi\phi_{h_C}(h)=\chi_C^{(i)}(h)$
for any $h\in Z_{u(C)}$. Consequently,
$$\chi_{\phi(C)}^{\prime(\phi_C(i))}=\chi_C^{(i)}(\phi\phi_{h_C})^{-1}, $$ i.e. $$\chi_{\phi(C)}^{\prime
(i)}=\chi_C^{(\phi_C^{-1}(i))}(\phi\phi_{h_C})^{-1}, $$ where
$C\in \mathcal {K}_r(G)$. That is,  (ii) holds.

Conversely,  if (i) and (ii) hold,  let
$$r'=\sum\limits_{C\in \mathcal {K}(G)}r_{\phi^{-1}(C)}C$$
and $$I_{\phi(C)}(r')=I_C(r)$$ for any $C\in\mathcal {K}(G)$,  $
\phi\in$Aut$G$. Obviously,  $\mathcal {K}_{r'}(G)=\{\phi(C) \mid
C\in \mathcal {K}_r(G)\}.$ For any $C\in \mathcal {K}_r(G), h_C\in
G,  \phi_{_C}\in S_{I_C(r)}=S_{r_{_C}}$,  let
$$u'(\phi(C))=\phi(h_Cu(C)h_C^{-1})$$ and
$$\chi_{\phi(C)}^{\prime
(i)}=\chi_C^{(\phi_C^{-1}(i))}(\phi\phi_{h_C})^{-1}$$ for any
$i\in I_{\phi(C)}(r').$ It is easy to check $RSC(G, r',
\overrightarrow{\chi'}, u')\cong RSC(G, r, \overrightarrow{\chi},
u)$.$\Box$

\begin{Corollary}\label{2.6}
If $G=S_n$ with $n\neq 6$,  then $RSC(G,  r,
\overrightarrow{\chi},  u)\cong RSC(G,  r',
\overrightarrow{\chi'},  u') $ if and only if

(i) $r'=r$;

(ii)
$\chi_C^{\prime(i)}=\chi_C^{(\phi_C^{-1}(i))}\phi^{-1}_{g_{_C}}$
and $u'(C)=\phi_{g_{_C}}(u(C))$ \ for any  $ C\in \mathcal
{K}_r(G)$,  $ i\in I_C(r) $,   $\phi\in \mathrm{Aut}G$ and
$g_{_C}\in G$.
\end{Corollary}

{\bf Proof.} By Proposition \ref{7.2}(ii),
$\mathrm{Aut}G=\mathrm{Inn}G$ and $\phi(C)=C$ for any $ \phi
\in$Aut$G,  C\in \mathcal {K}(G)$. There exists $g_{_C}\in G$ such
that  $\phi\phi_{h_C}=\phi_{g_{_C}}$ since $\phi\phi_{h_C}\in
\mathrm{Aut}G=\mathrm{Inn} G$ for any $h_C \in G,   \phi
\in$Aut$G,  C\in \mathcal {K}(G)$.  Considering Proposition
\ref{2.5},  we complete the proof.$\Box$

From now on,   we suppose that $G=S_n(n\neq 6)$. For a given
ramification $r$ of $G$,  let  $ \Omega (G,  r)$ be the set of all
RSC's of G with the ramification $r$,  namely,  $ \Omega (G,  r) :
= \{ (G,  r,  \overrightarrow{ \chi },  u) \mid (G,   r,
\overrightarrow{ \chi },  u) \mbox {\ is\ an  } RSC\}$. Let ${\cal
N}(G,  r)$ be the number of isomorphic classes in $ \Omega (G,
r)$. This article is mostly devoted to investigate the formula of
${\cal N}(G,  r)$.

\subsection{Action of Group }

Denote $RSC(G,  r,  \overrightarrow{\chi},  u)$ by
$\left(\left\{\chi^{(i)}_{C}\right\}_i,
u(C)\right)_{C\in\mathcal{K}_r(G)}$ for short. Let
$M:=\prod\limits_{C\in \mathcal {K}_r(G)}(\mbox{Aut}G\times
S_{I_C(r)})$ $=\prod\limits_{C\in \mathcal
{K}_r(G)}(\mbox{Inn}G\times S_{r_{_C}})$. Define the action of $M$
on $\Omega(G,  r)$ as follows:
\begin{eqnarray} \label {e2.1}
(\phi_{g_{_C}},
\phi_C)_C\circ\left(\left\{\chi^{(i)}_{C}\right\}_i, u(C)\right)_C
&=&
\left(\left\{\chi_C^{(\phi_C^{-1}(i))}\phi^{-1}_{g_{_C}}\right\}_i,
\phi_{g_{_C}}u(C) \right )_C
\end{eqnarray}
 In term of Corollary  \ref{2.6},  each orbit of $\Omega(G,
r)$ represents an isomorphic class of RSC's. As a result,  ${\cal
N}(G,  r)$ is equal to the number of orbits in $\Omega(G,  r )$.

We compute the the number of orbits by following steps.

\subsection{ The structure of centralizer in $S_n$ }

We recall the semidirect product and the wreath product of groups
(see \cite[Page 21]{char},  \cite [Page 268-272]{suzuki}, ).
\begin{Definition}\label{2.7}

Let $N$ and $K$ be groups, and assume that there exists a group
homomorphism $\alpha:K\rightarrow \mathrm{Aut}N$. The semidirect
product $N\rtimes_{\alpha} K$ of $N$ and $K$ with respect to
$\alpha$ is defined as follows:

(1) As set,  $N\rtimes_{\alpha} K$ is the Cartesian product of N
and K;

(2) The multiplication is given by
$$
(a, x)(b, y)=(a\alpha(x)(b), xy)
$$
for any a $\in$ N and x $\in$ N.
\end{Definition}

{\bf Remark.} Let $L=N\rtimes_{\alpha} K$. If we set
$$ \bar{N}=\{(a, 1) \mid a\in N\}, \bar{K}=\{(1, x) \mid x\in
K)\}, $$ then $ \bar{N}$ and $\bar{K}$ are subgroups of $L$ which
are isomorphic to $N$ and $K$,  respectively,  and will be
identified with $N$ and $K$.

\begin{Definition}\label{2.10}
 Let $A$ and $H$ be two  groups  and $H$ act on the set $X$. Assume that
 $B$ is the set  of all maps from $X$ to $A$. Define the multiplication of $B$ and
 the group homomorphism $\alpha $ from  $H$  to  Aut $B$ as follows:
 $$
(bb')(x)=b(x)b'(x) \hbox { \ and } \alpha (h)(b)(x)=b(h^{-1}\cdot
x)$$ for any $x\in X,  h\in H,  b,  b' \in B.$

 The  semidirect product
$B\rtimes_{\alpha} H$ is called the (general) wreath product of
$A$ and $H$,  written as  $A$\:wr\:$H$.

\end{Definition}

{\bf Remark.}\cite[Page 272]{suzuki} indicated that in Definition
\ref {2.10} $B\cong \prod\limits_{x\in X}A_x,  $ where $A_x= A$
for any $x\in X$. Hence $W=A$\:wr\:$H=B\rtimes
H=(\prod\limits_{x\in X}A_x)\rtimes H$. In particular. the wreath
product  plays  an important  role in the structure of centralizer
in $S_n$ when $H= S_n$ and $X = \{1,  2,  \cdots,  n\}$.

Now we keep on the work in \cite[Page 295-299 ]{suzuki}.
 Let $Y_i $ be the set of all letters which
belong to those cycles of length $i$ in the independent  cycle
decomposition of $\sigma$. Clearly,  $Y_{i}\bigcap
Y_{j}=\emptyset$ for $i\neq j$ and $\bigcup\limits_iY_i=\{1, 2,
\cdots, n \}$.

\begin{Lemma}\label{2.11}
 $\rho \in Z_{\sigma }$ if and only
if  $Y_{i}$ is $\rho$-invariant,  namely $\rho(Y_i)\subset Y_i$,
and
 the restriction $\rho_{i}$ of $\rho$ on $Y_{i}$ commutes with the
restriction $\sigma_{i}$ of $\sigma$ on $Y_{i}$ for  $i = 1,  2,
\cdots,  n$.
\end{Lemma}

{\bf Proof.}It follows from \cite[Page 295]{suzuki}.$\Box$

\begin{Proposition}\label{2.12}
If $\sigma$ and  $ \sigma_i$ is the same as above,  then

(i) $Z_{\sigma}= \prod\limits_{i}Z_{\sigma_i};$

(ii) $Z_{\sigma_i}\cong C_i$\:wr\:$S_{\lambda_i }\cong
\left(C_i\right)^{\lambda_i}\rtimes S_{\lambda_i} ,   $ where
$Z_{\sigma_i}$ is the centralizer of $\sigma_i$ in $S_{Y_i}$;
$\left(C_i\right)^{\lambda_i}=\overbrace{C_i\times \cdots\times
C_i}^{\lambda_i\mbox{个}}$; define $Z_{\sigma_i}=1$ and $
S_{\lambda_i}=1 $ when $\lambda_i=0$.
\end{Proposition}

{\bf Proof. }(i) It follows from Lemma  \ref {2.11}.

 (ii) Assume $$\tau=(a_{10}a_{11}\cdots a_{1,  l-1})(a_{20}\cdots a_{2,
l-1})\cdots (a_{m0}\cdots a_{m,  l-1}) $$ and
$$
X=\{1,  2,  \cdots,  m\},  \ \  Y=\{a_{ij} \mid  1 \le i \le m,  0
\le j \le l-1 \}.$$ Obviously,  in order to prove (ii),  it
suffices to show that $Z_{\tau}\cong C_l$\:wr\:$S_m$,  where
$Z_{\tau}$ is the centralizer of $\tau$ in $S_Y$. The second index
of  $a_{ij}$ ranges over $\{0, 1, \cdots, l-1\}$. It is convenient
to identify this set with $C_l=\mathbb{Z}/(l)$,  and to use
notation such as $a_{i, l+j}=a_{ij}$. We will define an
isomorphism $\varphi$ from $Z_{\tau}$ onto $C_l$\:wr\:$S_m$. For
any element $\rho$ of $Z_{\sigma}$,  we will define a permutation
$\theta(\rho)\in S_m$ and a function $f(\rho)$ from $X$ into $C_l$
by
$$\rho^{-1}(a_{i0})=a_{jk}, \qquad j=\theta(\rho)^{-1}(i), \qquad
f_{\rho}(i)=g^k, $$ where $g$ is the generator of $C_l$,  $1\le i
\le m$. Let
$$\varphi(\rho)=(f_{\rho},  \theta(\rho)). $$
Then $\varphi$ maps $Z_{\sigma}$ into $C_l$\:wr\:$S_m$. First we
 show that $\varphi$ is surjective. In fact,  for any $(f,
\theta) \in C_l\ \mathrm{wr}\ S_m$ and $i\in X$,  if $f (i) = g^k,
$ $\theta (j) = i$,  then there exists a permutation $\rho$ such
that $\rho ^{-1} (a_{ir})= a _{j,  k+r}$. It is easy to verify
that $\rho\in Z_{\sigma}$ and $\varphi (\rho ) = (f,  \theta)$.
Next we show that $\varphi$ is injective. If $\rho\in$
Ker$\varphi$,  then $\theta(\rho)=1$ and $f _{\rho}(i)=g^0$,  this
means that $\rho(a_{i0})=a_{i0}$ for all $i$. Hence,
$\rho(a_{ij})=a_{ij}$ for all $i$ and $j$ and $\rho=1_Y$,  which
implies that $\varphi$ is bijective. We complete the proof if
$\varphi$ is homomorphism,  which is proved below. For any $\rho,
\pi\in Z_{\tau}$,  suppose that
$$
(\rho\pi)^{-1}(a_{i0})=\pi^{-1}(a_{jk})=a_{s,  t+k},
$$
where $j=\theta(\rho)^{-1}(i), g^k=f_{\rho}(i),
s=\theta(\pi)^{-1}(j)=\theta(\pi)^{-1}\theta(\rho)^{-1}(i),
g^t=f_{\pi}(j)=f_{\pi}(\theta(\rho)^{-1}(i))$. Consequently,
$$\theta(\rho\pi)^{-1}(i)=s=(\theta(\rho)\theta(\pi))^{-1}(i),  \qquad
f_{\rho\pi}(i)=g^{t+k}=f_{\rho}(i)f_{\pi}(j)=f_{\rho}(i)f_{\pi}(\theta(\rho)^{-1}(i)).
$$

On the other hand,  in view of  Definition \ref {2.10},
$\theta(\rho)(f_{\pi})(i)=f_{\pi}(\theta(\rho)^{-1}(i))=f_{\rho\pi}(i)$,
whence
$$(f_{\rho},  \theta(\rho))(f_{\pi},
\theta(\pi))=(f_{\rho}\theta(\rho)(f_{\pi}),
\theta(\rho)\theta(\pi))=(f_{\rho\pi}, \theta(\rho\pi)).$$ It
follows that
$$\varphi(\rho\pi)=\varphi(\rho)\varphi(\pi), $$
which is just what we need.\ $\Box$

\begin{Corollary} \label{7.15}

{{ (i)}
$|Z_{\sigma}|=\lambda_1!\lambda_2!\cdots\lambda_n!1^{\lambda_1}2^{\lambda_2}\cdots
n^{\lambda_n}, $ where $\sigma$ is the same as above;

{(ii)} $\sum\limits_{h\in C}|Z_h|=|G|=n!$}  for any $ C\in
\mathcal {K}(G), $

\end{Corollary}

{\bf Proof. }  (i) It follows from \ref{2.12} that
$$
 \mid  Z_{\sigma} \mid  =\prod\limits_{i} \mid  Z_{\sigma_i} \mid
=\prod\limits_{i} \mid  \left(C_i\right)^{\lambda_i}\rtimes
S_{\lambda_i} \mid =\prod\limits_{i} \mid
\left(C_i\right)^{\lambda_i} \mid   \mid  S_{\lambda_i} \mid
=\prod\limits_{i}(i^{\lambda_i}\lambda_i!).
$$

(ii) Assume the type of $C$ is $1^{\lambda_1}2^{\lambda_2}\cdots
n^{\lambda_n}$. Combining (i) and Proposition \ref{7.2}(i),  we
have
$$
\sum\limits_{h\in C} \mid Z_h \mid
=(\lambda_1!\lambda_2!\cdots\lambda_n!1^{\lambda_1}2^{\lambda_2}\cdots
n^{\lambda_n})\frac{n!}{\lambda_1!\lambda_2!\cdots\lambda_n!1^{\lambda_1}2^{\lambda_2}\cdots
n^{\lambda_n}}=n! \ . \ \ \Box
$$

Proposition \ref{2.12} shows that $Z_{\sigma}$ is a direct product
of some wreath products  such as
 $C_l\:$wr\:$S_m$.
Denote by $((b_i)_i, \sigma)$ the element of
$C_l\:$wr\:$S_m=(C_l)^{m}\rtimes S_m$ where $b_i\in C_l, \sigma\in
S_m$. By Definition \ref{2.7} and \ref{2.10}, we have
\begin{eqnarray*}
((b_i)_i, \sigma)((b'_i)_i, \sigma')&=&((b_i)_i\cdot\sigma((b'_i)_i), \sigma\sigma')\nonumber\\
&=&((b_i)_i\cdot(b'_{\sigma^{-1}(i)})_i, \sigma\sigma')\nonumber\\
&=&((b_ib'_{\sigma^{-1}(i)})_i, \sigma\sigma')
\end{eqnarray*}
and
\begin{equation*}
((b_i)_i, \sigma)^{-1}=((b^{-1}_{\sigma(i)})_i, \sigma^{-1}).
\end{equation*}

\subsection { Characters of $Z_{\sigma}$}

\begin{Lemma}\label{2.14} If
$H=H_1\times\cdots\times H_r$,  then

(i) $H'=H'_1\times\cdots\times H'_r$;

(ii) $H/H'\cong \prod\limits_{i=1}^{r}H_i/H'_i.$
\end{Lemma}
{\bf Proof.} (i) It follows from  \cite  [ Page 77] {Zh82};

(ii) It follows from \cite[Corollary 8.11]{Hu74} and (i). \ $\Box$

According to Proposition \ref{2.12}(i) and Lemma \ref{2.14}(ii),
to investigate  the structure of $Z_{\sigma_i}/Z'_{\sigma_i}$,  we
have to study $(C_l\:$wr\:$S_m)/(C_l\:$wr\:$S_m)'$.

\begin{Lemma}\label{2.15}
Let $B=(C_l)^{m}$,   $\bar{B}=\left\{(b_i)_i\in B \mid
\prod\limits_{i=1}^mb_i=1\right\} $ and
$W=C_l\:$wr\:$S_m=B\rtimes S_m$. Then $W'=\bar{B}\rtimes
S'_m=\bar{B}\rtimes A_m$.
\end{Lemma}

{\bf Proof. } Obviously,  $\bar{B}\lhd B$.

For any  $((b_i)_i,  \sigma),  ((b'_i)_i,  \sigma')\in W, $ see
\begin{eqnarray*}
& &((b_i)_i,  \sigma)^{-1}((b'_i)_i,  \sigma')^{-1}((b_i)_i,  \sigma)((b'_i)_i,  \sigma')\\
& = &((b^{-1}_{\sigma(i)})_i,  \sigma^{-1})((b^{\prime-1}_{\sigma'(i)})_i,  \sigma^{\prime-1})((b_ib'_{\sigma^{-1}(i)})_i,  \sigma\sigma')\\
& = &((b^{-1}_{\sigma(i)}b^{\prime-1}_{\sigma'\sigma(i)})_i,  \sigma^{-1}\sigma^{\prime -1})((b_ib'_{\sigma^{-1}(i)})_i,  \sigma\sigma')\\
& =
&((b^{-1}_{\sigma(i)}b^{\prime-1}_{\sigma'\sigma(i)}b_{\sigma'\sigma(i)}b'_{\sigma^{-1}\sigma'\sigma(i)})_i,
\sigma^{-1}\sigma^{\prime-1}\sigma\sigma').
\end{eqnarray*}
Considering $\prod\limits_ib_{\tau(i)}=\prod\limits_ib_i$
 for any $ \tau \in S_m$,  we have
\begin{eqnarray*}
&&\prod\limits_i\left(b^{-1}_{\sigma(i)}b^{\prime-1}_{\sigma'\sigma(i)}b_{\sigma'\sigma(i)}b'_{\sigma^{-1}\sigma'\sigma(i)}\right)\\
&=&\prod\limits_ib^{-1}_i\prod\limits_ib^{\prime-1}_i\prod\limits_ib_i\prod\limits_ib'_i\\
&=&1
\end{eqnarray*}
and
$b^{-1}_{\sigma(i)}b^{\prime-1}_{\sigma'\sigma(i)}b_{\sigma'\sigma(i)}b'_{\sigma^{-1}\sigma'\sigma(i)}\in
\bar{B}. $  By Proposition  \ref {7.2}(iii),
$\sigma^{-1}\sigma^{\prime-1}\sigma\sigma'\in S'_m=A_m. $ Thus
$$((b_i)_i,  \sigma)^{-1}((b'_i)_i,  \sigma')^{-1}((b_i)_i,  \sigma)((b'_i)_i,  \sigma')\in \bar{B}\rtimes
A_m,  $$ which implies $W'\subset \bar{B}\rtimes A_m.$

Conversely,  for any  $ b=((b_i)_i,  1)\in \bar{B}$,  then
$\prod\limits_{i=1}^mb_i=1$. Let $\sigma=(1,  2,  \cdots,  m),
b'_1=1,  b'_i=\prod\limits_{j=1}^{i-1}b_j\;(2\leq i\leq m)$,   See
\begin{eqnarray*}
((b'_i)_i,  1)^{-1}(1,  \sigma)^{-1}((b'_i)_i,  1)(1,  \sigma)
&=&((b^{\prime-1}_i)_i,  1)(1,  \sigma^{-1})((b'_i)_i,  \sigma)\\
&=&((b^{\prime-1}_i,  b'_{\sigma(i)})_i,  1)\\
&=&((\prod\limits_{j=1}^{i-1}b^{-1}_j\prod\limits_{j=1}^{i}b_j)_i,  1)\\
&=&((b_i)_i,  1)=b.
\end{eqnarray*}
Therefore,   $b$ is a commutator of $((b'_i)_i,  1)$ and $(1,
\sigma)$,  i.e. $b\in W'.$ This show $\bar{B}\subset W'.$

Furthermore,  for any  $ \sigma,  \sigma'\in S_m,  $,  we have
$$(1,  \sigma)^{-1}(1,  \sigma')^{-1}(1,  \sigma)(1,  \sigma')=(1,  \sigma^{-1}\sigma^{\prime-1}\sigma\sigma'), $$
which implies $(1,  \sigma^{-1}\sigma^{\prime-1}\sigma\sigma')\in
W'$ and $S'_m=A_m\subset W'.$

Consequently,   $\bar{B}\rtimes A_m \subset W'. $$\Box$

\begin{Proposition}\label{2.17}
 (i) Let $B$ and  $\bar{B}$ be the same as in Lemma  \ref {2.15}. Then $B/\bar{B}\cong C_l;$

(ii) Let $W$ and  $W'$ be the same as in Lemma \ref {2.15}. Then $$
W/W'\cong\left\{
\begin{array}{lll}
C_l\times C_2& m\geq2\\
C_l &m =1\\
1 & m =0
\end{array}
\right. ; $$

(iii) Let $\sigma$ be the same as in Proposition \ref{2.12}. Then
$$Z_{\sigma}/Z'_{\sigma}\cong\left(\prod\limits_{\lambda_i=1}C_{i}\right)
\left(\prod\limits_{\lambda_i\geq2}(C_{i}\times C_2)\right).$$
\end{Proposition}

{\bf Proof. } (i) It is clear that  $\nu : B\rightarrow C_l$ by
sending $(b_i)_i$ to $ \prod\limits_ib_i$  is a group homomorphism
with Ker$\nu=\bar{B}$. Thus (i) holds.

(ii) By (i) and Proposition \ref{7.2}(iii),  we have to show
$W/W'\cong (B/\bar{B})\times (S_{m}/A_{m}). $ Define
\begin{align*}
\psi:W &\rightarrow (B/\bar{B})\times (S_{m}/A_{m})\\
(b,  \sigma) &\mapsto (\bar{b},  \bar{\sigma}),
\end{align*}
where $\bar{b}$ and $\bar{\sigma}$ are the images of $b$ and
$\sigma$ under canonical epimorphisms $B\rightarrow B/\bar{B}$ and
$S_m\rightarrow S_{m}/A_{m}$,  respectively.  Therefore $$
\psi((b',  \sigma')(b,  \sigma))=(\overline{b'\sigma'(b)},
\overline{\sigma'\sigma})=(\bar{b'}\overline{\sigma'(b)},
\bar{\sigma'}\bar{\sigma}). $$  Let $b=(b_i)_i$. Since
$\nu(\sigma'(b))=\nu((b_{\sigma^{\prime-1}(i)})_i)=\prod\limits_ib_{\sigma^{\prime-1}(i)}=\prod\limits_ib_i=\nu(b)$,
we have $\overline{\sigma'(b)}=\bar{b}$. Consequently,
$$\psi((\sigma',  b')(\sigma,  b))=(\bar{b'}\overline{\sigma(b)},  \bar{\sigma'}\bar{\sigma})=(\bar{b'}\bar{b},  \bar{\sigma'}
\bar{\sigma},  )=\psi((b',  \sigma'))\psi((b,  \sigma)). $$ This
show that $\psi$ is a group homomorphism. Obviously,   it is
surjective and Ker$\psi=\bar{B}\rtimes A_{m}=W'$. we complete the
proof of (ii).

(iii) By Proposition \ref{2.12}(i) and Lemma \ref{2.14}(ii),
$Z_{\sigma}/Z'_{\sigma}=(\prod\limits_i
Z_{\sigma_i})/(\prod\limits_i Z'_{\sigma_i})\cong\prod\limits_i
(Z_{\sigma_i}/Z'_{\sigma_i}).$ Applying (ii),  we complete the
proof of (iii). $\Box$

\begin{Corollary}\label{2.18}
If $\sigma \in C$ is a permutation  of type
$1^{\lambda_1}2^{\lambda_2}\cdots n^{\lambda_n}$,  then $ \mid \!
\widehat {Z_{\sigma }}\! \mid
=\left(\prod\limits_{\lambda_i=1}i\right)\left(\prod\limits_{\lambda_i\geq2}2i\right)$,
written as  $\gamma _C$.
\end{Corollary}

{\bf Proof. } It follows from Corollary \ref{7.10} and Proposition
\ref{2.17}(iii). $\Box$

{\bf Remark.} Obviously,  $\gamma_{_C}$ only depends on the
conjugate class $C$．

\subsection { Fixed Point Set $F_{(\phi_{g_{_C}},  \phi_C)_C}$  }
  Let $(\phi_{g_{_C}}, \phi_{_C})_C\in M$. If  $\left(\left\{\chi^{(i)}_{C}\right\}_i,  u(C)\right)_C \in
F_{(\phi_{g_{_C}},  \phi_C)_C},  $ then,  according to (\ref
{e2.1}),  we have
\begin{eqnarray*}
(\phi_{g_{_C}},
\phi_C)_C\circ\left(\left\{\chi^{(i)}_{C}\right\}_i,
u(C)\right)_C&=&\left(\left\{\chi_C^{(\phi_C^{-1}(i))}\phi^{-1}_{g_{_C}}\right\}_i,
\phi_{g_{_C}}u(C)
\right)_C\\
&=&\left(\left\{\chi^{(i)}_{C}\right\}_i,  u(C)\right)_C.
\end{eqnarray*}
Consequently,
\begin{eqnarray} \label {e2.2}
\phi_{g_{_C}}(u(C))=u(C)
\end{eqnarray}
and
\begin{eqnarray} \label {e2.3}
\chi_C^{(\phi_C^{-1}(i))}\phi^{-1}_{g_{_C}}=\chi^{(i)}_{C}
\end{eqnarray}
for any $ C \in\mathcal {K}_r(G)$,  $i\in I_C(r)$.

 (\ref
{e2.2}) implies $u(C)\in Z_{g_{_C}}$ and $g_{_C}\in Z_{u(C)}$.
Considering $u(C)\in C$,  we have $u(C)\in Z_{g_{_C}}\cap C$.
Conversely,  if  $\sigma\in Z_{g_{_C}}\cap C$,  then $u(C)=\sigma$
satisfies (\ref {e2.2}). Consequently,  the number of $u(C)$'s
which satisfy (\ref {e2.2}) is  $\mid Z_{g_{_C}}\cap C \mid, $
written
$$\beta_{g_{_C}}:= \mid Z_{g_{_C}}\cap C \mid .$$

(\ref {e2.3}) is equivalent to
\begin{eqnarray} \label {e2.4}
\chi_C^{(i)}\phi^{-1}_{g_{_C}}=\chi^{(\phi_C(i))}_{C},
\end{eqnarray}
where $\phi_C\in S_{I_C(r)}=S_{r_{_C}}$. It is clear that (\ref
{e2.4}) holds if and only if the following (\ref {e2.4'}) holds
for every  cycle,  such as  $\tau=(i_1,  \cdots,  i_r)$,   of
$\sigma$:
\begin{equation}\tag{\ref{e2.4}$'$}\label{e2.4'}
\chi_C^{(i)}\phi^{-1}_{g_{_C}}=\chi^{(\tau(i))}_{C}, \ i=i_1,
\ldots, i_r.
\end{equation}

By (\ref{e2.4'}),  we have
\begin{eqnarray}\label {e2.5}
\chi^{(i_2)}_{C}&=&\chi_C^{(i_1)}\phi^{-1}_{g_{_C}},  \nonumber\\
\chi^{(i_3)}_{C}&=&\chi_C^{(i_2)}\phi^{-1}_{g_{_C}}=\chi_C^{(i_1)}(\phi^{-1}_{g_{_C}})^2,  \nonumber\\
&\cdots&\nonumber\\
\chi^{(i_r)}_{C}&=&\chi_C^{(i_{r-1})}\phi^{-1}_{g_{_C}}=\chi_C^{(i_1)}(\phi^{-1}_{g_{_C}})^{r-1},  \nonumber\\
\chi^{(i_1)}_{C}&=&\chi_C^{(i_{r})}\phi^{-1}_{g_{_C}}=\chi_C^{(i_1)}(\phi^{-1}_{g_{_C}})^{r}.
\end{eqnarray}

We can obtain $\chi^{(i_2)}_{C},  \cdots,  \chi^{(i_r)}_{C}$ when
we obtain $\chi_C^{(i_1)}$. Notice
$(\phi^{-1}_{g_{_C}})^{r}=(\phi^{r}_{g_{_C}})^{-1}=(\phi_{g^r_{_C}})^{-1}$.
Therefore  (\ref {e2.5}) is equivalent to $$
\chi^{(i_1)}_{C}\phi_{g^r_{_C}}=\chi^{(i_1)}_{C},  $$ i.e.
\begin{eqnarray} \label{e2.6}
\chi^{(i_1)}_{C}\left(g^r_{_C}h(g^r_{_C})^{-1}\right)=\chi^{(i_1)}_{C}(h)
\end{eqnarray} for any $
h\in Z_{u(C)}$.  Since $g_{_C}\in Z_{u(C)}$ implies $g_{_C}^r\in
Z_{u(C)}$,  we have $$
\chi^{(i_1)}_{C}\left(g^r_{_C}h(g^r_{_C})^{-1}\right)=
\chi^{(i_1)}_{C}\left(g^r_{_C}\right)\chi^{(i_1)}_{C}\left(h\right)\chi^{(i_1)}_{C}\left((g^r_{_C})^{-1}\right)
=\chi^{(i_1)}_{C}\left(h\right).
$$
This shows (\ref {e2.6}) holds for any characters in $Z_{u(C)}$,
By Corollary \ref{2.18},  there exist $\gamma_{_C}$ characters
$\chi_C^{(i_1)}$'s such that  (\ref {e2.5}) holds. Assume that
$\phi_C$  is  written as multiplication of $k_{\phi_C}$
independent cycles. Thus given a $u(C)$,  there exist
$\gamma_C^{k_{\phi_C}}$ elements $\{\chi^{(i)}_C\}_{i}$'s such
that  (\ref {e2.4}) holds. Consequently,  for any $C\in \mathcal
{K}_r(G)$,  there exist $\gamma_C^{k_{\phi_C}}\beta_{g_{_C}}$
distinct elements $\left(\left\{\chi^{(i)}_{C}\right\}_i,
u(C)\right)$'s  such that both (\ref {e2.2})and (\ref {e2.3})
hold.

 Consequently,  we have
\begin{equation}\label{e2.77}
  \mid  F_{(\phi_{g_{_C}},
\phi_C)_C} \mid  =\prod\limits_{ C\in \mathcal
{K}_r(G)}\left(\gamma_C^{k_{\phi_C}}\beta_{g_{_C}}\right)
\end{equation}
for any  $\phi_{g_{_C}},  \phi_C)_C\in M.$
\subsection{The Formula Computing the Number of Isomorphic Classes}

Applying Burnside's Lemma and (\ref{e2.77}),  we have
\begin{eqnarray}\label {e2.7} {\cal N} (G, r)
&=&\frac{1}{ \mid  M \mid  }\sum\limits_{(\phi_{g_{_C}},
\phi_C)_C\in
M} \mid  F_{(\phi_{g_{_C}},  \phi_C)_C} \mid  \nonumber\\
&=&\frac{1}{ \mid  M \mid  }\sum\limits_{(\phi_{g_{_C}},
\phi_C)_C\in M}\prod\limits_{ C\in \mathcal
{K}_r(G)}\left(\gamma_C^{k_{\phi_C}}\beta_{g_{_C}}\right)\nonumber\\
&=&\frac{1}{ \mid  M \mid  }\prod\limits_{ C\in \mathcal
{K}_r(G)}\sum\limits_{\phi_{g_{_C}}\in \mathrm{Inn}G\atop \phi_C\in
S_{r_C}}\left(\gamma_C^{k_{\phi_C}}\beta_{g_{_C}}\right)\nonumber\\
&=&\frac{1}{ \mid  M \mid  }\prod\limits_{ C\in \mathcal
{K}_r(G)}\sum\limits_{\phi_{g_{_C}}\in
\mathrm{Inn}G}\sum\limits_{\phi_C\in
S_{r_C}}\left(\gamma_C^{k_{\phi_C}}\beta_{g_{_C}}\right)\nonumber\\
&=&\frac{1}{ \mid  M \mid  }\prod\limits_{ C\in \mathcal
{K}_r(G)}\left(\sum\limits_{\phi_{g_{_C}}\in
\mathrm{Inn}G}\beta_{g_{_C}}\sum\limits_{\phi_C\in
S_{r_C}}\gamma_C^{k_{\phi_C}}\right)
\end{eqnarray}
where  $\beta_{g_{_C}}$ and $k_{\phi_C}$ are the same as above.  $
\mid  M \mid =\prod\limits_{C\in \mathcal {K}_r(G)}(n!r_{_C}!)$
when $n \not=2, 6$;
 $ \mid  M \mid  =\prod\limits_{C\in \mathcal {K}_r(G)}(r_{_C}!)$
 when $n=2$.

Next we simplify formula (\ref {e2.7}). In  (\ref {e2.7}),  since
$\gamma_C$ depends only on conjugate class,  it can be computed by
means of Corollary \ref{2.18}. Notice that $k_{\phi_C}$ denotes
the number of independent cycles in independent cycle
decomposition
 of $\phi _C$ in $S_{r_{_C}}$.

\begin{Definition}\label{2.19} (See \cite[pages 292-295]{Ri} )
Let
$$
[x]_n=x(x-1)(x-2)\cdots (x-n+1)\qquad n=1,  2,  \cdots.
$$ Obviously,  $[x]_n$ is a polynomial  with degree $n$. Let  $(-1)^{n-k}s(n,  k)$ denote the coefficient
of $x^k$ in $[x]_n$ and  $s(n,  k)$ is called the 1st Stirling
number. That is,
$$
[x]_n=\sum\limits_{k=1}^n(-1)^{n-k}s(n,  k)x^k.
$$
\end{Definition}

\cite[Theorem 8.2.9]{Ri} obtained the meaning about the 1st Stirling
number:
\begin{Lemma}\label{2.20} (\cite[Theorem 8.2.9]{Ri})
There exactly exist $s(n,  k)$ permutations in $S_n$,   of  which
the numbers of independent cycles in independent cycle
decomposition
 are $k$. That is,  $s(n,  k)$ = $\mid \! \{ \tau \in S_n \mid$ the number of  independent cycles in independent cycle
 decomposition of $\tau$ is $k \}\! \mid$. \end{Lemma}

About $\beta_{g}$,  we have

\begin{Lemma}\label{2.21}
Let $G=S_n$.  Then $\sum\limits_{g\in G}\beta_{g}=\sum\limits_{g\in
G} \mid  Z_g\cap C \mid  = \mid G \mid =n!$  for any $C\in\mathcal
{K}(G)$.
\end{Lemma}

{\bf Proof. } For any  $C\in \mathcal {K}(G),  $ we have
\begin{eqnarray*}
\sum\limits_{g\in G} \mid  Z_g\cap C \mid  &=& \sum\limits_{g\in
G}\left |  \bigcup\limits_{h\in C}\left(Z_g\cap
\{h\}\right)\right | \\
&=&\sum\limits_{g\in G}\sum\limits_{h\in
C}\left | Z_g\cap\{h\}\right |  \\
&=&\sum\limits_{h\in C}\sum\limits_{g\in
G}\left | Z_g\cap\{h\}\right | \\
&=&\sum\limits_{h\in C} \mid  Z_h \mid  \\
&=&n!  \ \ ( \mbox { \hbox { by Corollary }\ref {7.15}(ii)} ). \ \ \
\Box
\end{eqnarray*}

\begin{Theorem}\label{2.22}
Let $G=S_n$ with $n\neq6$. Then
\begin{eqnarray} \label {e2.8}
{\cal N} (G,  r)=\prod\limits_{C\in \mathcal
{K}_r(G)}{\gamma_C+r_C-1 \choose r_C}. \end  {eqnarray}
\end{Theorem}

{\bf Proof. } We shall show in following  two cases.

(i)  $n \not= 2,  6$. By Proposition \ref{7.2}(ii),  Inn$G=G$.
Applying (\ref{e2.7}),  we have
\begin{eqnarray*}
{\cal N} (S_n,   r) &=&\frac{1}{\prod\limits_{ C\in \mathcal
{K}_r(G)}(n!r_C!)}\prod\limits_{ C\in \mathcal
{K}_r(G)}\left(\sum\limits_{g_{_C}\in
G}\beta_{g_{_C}}\sum\limits_{\phi_C\in
S_{r_C}}\gamma_C^{k_{\phi_C}}\right)\\
&=&\prod\limits_{ C\in
\mathcal {K}_r(G)}\left(\frac{1}{n!}\sum\limits_{g_{_C}\in
G}\beta_{g_{_C}}\right)\prod\limits_{ C\in \mathcal
{K}_r(G)}\left(\frac{1}{r_C!}\sum\limits_{\phi_C\in
S_{r_C}}\gamma_C^{k_{\phi_C}}\right)\\
&=&\prod\limits_{ C\in \mathcal
{K}_r(G)}\left(\frac{1}{n!}n!\right)\prod\limits_{ C\in \mathcal
{K}_r(G)}\left(\frac{1}{r_C!}\sum\limits_{k=1}^{r_{_C}}s(r_{_C},
k)\gamma_C^k\right)
\ \ ( \hbox {by Lemma  \ref {2.21},   \ref {2.20}} ) \\
&=&\prod\limits_{ C\in \mathcal
{K}_r(G)}\left(\frac{1}{r_C!}(-1)^{r_{_C}}
\sum\limits_{k=1}^{r_{_C}}(-1)^{r_{_C}-k}s(r_{_C},  k)(-\gamma_C)^k\right)\ \ \\
&=&\prod\limits_{ C\in \mathcal
{K}_r(G)}\left(\frac{1}{r_C!}(-1)^{r_{_C}}(-\gamma_C)(-\gamma_C-1)\cdots(-\gamma_C-r_{_C}+1)\right)(\hbox{by Definition \ref{2.19}})\\
&=&\prod\limits_{ C\in \mathcal
{K}_r(G)}\frac{(\gamma_C+r_{_C}-1)\cdots(\gamma_C+1)\gamma_C}{r_{_C}!}\\
&=&\prod\limits_{C\in \mathcal {K}_r(G)}{\gamma_C+r_C-1 \choose
r_C}. \ \Box
\end{eqnarray*}

(ii) $n =2$. In this case we have  Inn$G=\{1\}$,
$\phi_{g_{_C}}=1$ and  $\beta_{g_{_C}}=1$ in (\ref{e2.7}). Thus
\begin{eqnarray*}
{\cal N} (S_n,   r) &=&\frac{1}{\prod\limits_{ C\in \mathcal
{K}_r(G)}(r_C!)}\prod\limits_{ C\in \mathcal
{K}_r(G)}\sum\limits_{\phi_C\in
S_{r_C}}\left(\gamma_C^{k_{\phi_C}}\right)\\
&=&\prod\limits_{ C\in \mathcal
{K}_r(G)}\left(\frac{1}{r_C!}\sum\limits_{\phi_C\in
S_{r_C}}\gamma_C^{k_{\phi_C}}\right)\ \ (\hbox{ similar to the computation of (i)})\\
&=&\prod\limits_{C\in \mathcal {K}_r(G)}{\gamma_C+r_C-1 \choose
r_C}. \ \Box
\end{eqnarray*}

\begin{Corollary}\label{2.23}
If $r=\sum\limits_{C\in \mathcal {K}(G)}C$,  i.e. $r_{_C}=1 $ for
any $ C\in \mathcal {K}(G)$,  then
\begin{eqnarray} \label {e2.9}
{\cal N} (G,  r)=\prod\limits_{C\in \mathcal {K}(G)}\gamma_C.
\end{eqnarray}
\end{Corollary}

{\bf Proof. } Considering $\mathcal {K}_r(G)=\mathcal {K}(G)$ and
applying Theorem \ref{2.22},  we can complete the proof. $\Box$

\begin{Corollary}\label{2.24}
If $C_0=\{(1)\}$ and $\mathcal {K}_r(G)=\{C_0\}, $ then
$$ {\cal N} (G,  r)=\left\{
\begin{array}{lll}
1& n=1\\
r_{_{C_0}}+1 & n\geq2
\end{array}
\right..$$
\end{Corollary}

{\bf Proof. } By Corollary \ref{2.18},  we have $\gamma_{_{C_0}}=1
$ when $n=1$;  $\gamma_{_{C_0}}=2 $ when $n\geq2$. Applying
Theorem \ref{2.22},  we can complete the proof. $\Box$

\subsection{Represetative}
 We have obtained  the number of isomorphic classes,  now  we give the
 representatives

\begin{Definition}\label{7.23}
Given an $RSC (G,   r,   \overrightarrow{\chi },   u)$,  written
$\widehat {Z_{u(C)}}=\{ \xi _{u(C)}^{(i)} \mid i = 1,  2,  \cdots,
\gamma _C\}$  and \ $n_i : =$ $ \mid \{ j \mid \chi _C ^{(j)} =
\xi _{u(C)}^{(i)} \} \! \mid, $ \ for any $C\in \mathcal {K}_r(G)$
\ \ and \ \ \ $1\le i \le \gamma _{u(C)}\  $ ,  then    $\{(n_1,
n_2,  \cdots,  n_{\gamma _C})\}_{C\in\mathcal {K}_r(G)}$ is called
the type of $RSC (G,  r,  \overrightarrow{\chi },  $ $ u)$.
\end{Definition}

\begin {Lemma} \label {2.30}
 $RSC (G,   r,   \overrightarrow {\chi},   u)$ and  $RSC (G,   r,
\overrightarrow {\chi'},   u)$ are isomorphic if and only if they
have the same type.

\end {Lemma}
{\bf Proof. } If  $RSC (G,   r,   \overrightarrow {\chi},   u)$
and $RSC (G,   r,   \overrightarrow {\chi'},   u)$ have the same
type,  then there exists a bijective map  $\phi _C \in S_{r_C}$
such that
 $\chi'{} ^{(\phi _C(i))}_C = \chi
^{(i)}_C$ for any  $i \in I_C(r)$,   $C \in {\cal K}_r (G)$. Let
$\phi = id _G$ and   $g_{_C} = 1$. It follows from Corollary
\ref{2.6} that they are isomorphic.

Conversely,  if  $RSC (G,   r,   \overrightarrow {\chi},   u)$ and
$RSC (G,  r,  \overrightarrow {\chi'},   u)$ are isomorphic, then,
by Corollary \ref{2.6},   there exist
 $\phi \in \mathrm{Aut }G$ and    $\phi _C \in S_{r_C}$,  as well as
 $g_{_C}\in G$ such that  $u(C)=\phi_{g_{_C}}(u(C))$ and
$\chi_{C}^{\prime(\phi_C(i))}\phi_{g_{_C}}(h)=\chi_C^{(i)}(h)$ for
any
 $C \in {\cal K}_r(G)$,   $h \in Z_{u(C)}$,    $i\in I_C(r)$.
Thus
  $g_C \in  Z_{u(C)}$ and  $\chi_{C}^{\prime(\phi_C(i))}\phi_{g_{_C}}(h)$ $=\chi '{}_{C}^{(\phi_C(i))}(g_C h g_C
^{-1})$ $=\chi_{C}^{\prime(\phi_C(i))}(h)=\chi_C^{(i)}(h) $ for
any $h \in Z_{u(C)}$,  i.e.
$\chi_{C}^{\prime(\phi_C(i))}=\chi_C^{(i)}.$ This implies that
$RSC (G,   r,  \overrightarrow {\chi},  u)$  and $RSC (G,  r,
\overrightarrow {\chi'},  u)$ have the same  type. $\Box$

\begin{Theorem}\label{2.31}
Let $G=S_n(n\neq6)$ and  $u_0$ be a map from ${\cal K}(G)
\rightarrow G$ with  $u_0(C) \in C$  for any $ C\in {\cal K}(G)$.
Let $\bar {\Omega }(G,   r,   u_0)$ denote the set consisting of
all elements with  distinct type in  $\{ (G,  r,
\overrightarrow{\chi },  u_0) \mid (G,  r,  \overrightarrow{\chi
},  u_0)$ is an  RSC  $\}$. Then $\bar {\Omega }(G,
 r,   u_0)$ becomes the representative set of  $\Omega (G,   r)$.
\end{Theorem}

{\bf Proof. } According to Theorem \ref {2.22} and \cite  {cao},
the number of elements in  $\bar {\Omega }(G,   r,   u_0)$  is the
same as the number of isomorphic classes in  $\Omega (G,  r)$.
Applying Lemma \ref {2.30} we complete the proof. $\Box$

For example,  applying Corollary \ref{2.23} we compute ${\cal
N}(G,  r)$ for $ G=S_3,  S_4,  S_5$ and $r=\sum\limits_{C\in
\mathcal {K}(G)}C$.

\begin{Example}\label{7.16}
(i) Let  $G = S_3$ and $r=\sum\limits_{C\in \mathcal {K}(G)}C$. Then
\begin{center}
\begin{tabular}{c | c | c |  c}
type  &$1^3$ & $1^12^1$ & $3^1$\\
\hline
 $\gamma_C$ & 2 & 2 & 3\\
\end{tabular}
\end{center}

By (\ref {e2.9}),  ${\cal N}(S_3,  r)=12.$ We explicitly   write
the type of representative elements as follows:\\$\{(1,  0)_{1^3},
(1,  0)_{1^12^1},  (1,  0,  0)_{3^1}\},  \{(0,  1)_{1^3},  (1,
0)_{1^12^1},  (1,  0,  0)_{3^1}\},  \{(1,  0)_{1^3},  (0,
1)_{1^12^1},  (1,  0,  0)_{3^1}\},
\\\{(0,  1)_{1^3},  (0,  1)_{1^12^1},  (1,  0,  0)_{3^1}\},
\{(1,  0)_{1^3},  (1,  0)_{1^12^1},  (0,  1,  0)_{3^1}\},  \{(1,  0)_{1^3},  (1,  0)_{1^12^1},  (0,  0,  1)_{3^1}\},  \\
\{(0,  1)_{1^3},  (1,  0)_{1^12^1},  (0,  1,  0)_{3^1}\},  \{(0,
1)_{1^3},  (1,  0)_{1^12^1},  (0,  0,  1)_{3^1}\},  \{(1,
0)_{1^3},  (0,  1)_{1^12^1},  (0,  1,  0)_{3^1}\},  \\\{(1,
0)_{1^3},  (0,  1)_{1^12^1},  (0,  0,  1)_{3^1}\},  \{(0,
1)_{1^3},  (0,  1)_{1^12^1},  (0,  1,  0)_{3^1}\},  \{(0,
1)_{1^3},  (0,  1)_{1^12^1},  (0,  0,  1)_{3^1}\}. $

(ii)  Let  $G = S_4$ and $r=\sum\limits_{C\in \mathcal {K}(G)}C$.
Then

\begin{center}
\begin{tabular}{c |  c |  c |  c |  c |  c}
type  & $1^4$ & $1^13^1$ & $1^22^1$&$2^2$&$4^1$\\
\hline
 $\gamma_C$ & 2 & 3 & 4& 4&4\\
\end{tabular}
\end{center}
By (\ref {e2.9}),  ${\cal N}(S_4,   r)=384. $

(iii)  Let  $G = S_4$ and $r=\sum\limits_{C\in \mathcal {K}(G)}C$.
Then
\begin{center}
\begin{tabular}{c | c |  c |  c | c | c |  c | c}
type  & $1^5$ & $1^14^1$ & $1^23^1$&$1^32^1$&$1^12^2$&$2^13^1$&$5^1$\\
\hline
 $\gamma_C$ & 2 & 4& 6& 4&4&6&5\\
\end{tabular}
\end{center}
By (\ref {e2.9}),   ${\cal N}(S_5,   r)=23040. $
\end{Example}

\end{document}